\documentclass{article}%
\usepackage{graphicx}
\usepackage{amsmath}
\usepackage{amsfonts}
\usepackage{amssymb}%
\providecommand{\U}[1]{\protect\rule{.1in}{.1in}}
\providecommand{\U}[1]{\protect\rule{.1in}{.1in}}
\newtheorem{theorem}{Theorem}

\newtheorem{lemma}[theorem]{Lemma}

\newtheorem{proposition}[theorem]{Proposition}
\newtheorem{remark}[theorem]{Remark}

\newenvironment{proof}[1][Proof]{\textbf{#1.} }{\ \rule{0.5em}{0.5em}}
\begin{document}

\title{On a new multivariate sampling paradigm and a polyspline Shannon function}
\author{Ognyan Kounchev and Hermann Render}
\maketitle

\begin{abstract}
In \cite{okbook} and \cite{kounchevrenderJAT} we have introduced and studied a
new paradigm for cardinal interpolation which is related to the theory of
multivariate polysplines. In the present paper we show that this is related to
a new sampling paradigm in the multivariate case, whereas we obtain a Shannon
type function $S\left(  x\right)  $ and the following Shannon type formula:
\[
f\left(  r\theta\right)  =\sum_{j=-\infty}^{\infty}\int_{\mathbb{S}^{n-1}%
}S\left(  e^{-j}r\theta\right)  f\left(  e^{j}\theta\right)  d\theta.
\]
This formula relies upon infinitely many Shannon type formulas for the
exponential splines arising from the radial part of the polyharmonic operator
$\Delta^{p}$ for fixed $p\geq1$.

Acknowledgement. The first and the second author have been partially supported
by the Institutes partnership project with the Alexander von Humboldt
Foundation. The first has been partially sponsored by the Greek-Bulgarian
bilateral project BGr-17, and the second author by Grant MTM2006-13000-C03-03
of the D.G.I. of Spain.

\end{abstract}

\section{Introduction}

The classical Shannon-Kotelnikov-Whittaker formula represents a function by
means of the following series
\[
f\left(  t\right)  =\sum_{j=-\infty}^{\infty}f\left(  jT\right)  \frac
{\sin\sigma\left(  t-jT\right)  }{\sigma\left(  t-jT\right)  }%
\]
where $T=\frac{\pi}{\sigma}$ and $f$ is a continuous function in $L_{2}\left(
\mathbb{R}\right)  $ having a Fourier transform with support in $\left[
-\pi,\pi\right]  .$

There have been many generalizations of this formula, see \cite{higgins},
\cite{ciaurri}, and references therein. There exists intimate relation between
sampling theorems and Wavelet Analysis which have been studied exhaustively in
\cite{walter}, see also the references therein.

In the present paper we will consider a multivariate generalization of the
Shannon theory which is based on a \textbf{semi-orthogonal} Wavelet Analysis
using polysplines, and which has been recently developed in \cite{okbook},
\cite{kounchevrenderappr}, \cite{aldazkounchevrender},
\cite{kounchevrenderJAT}, \cite{kounchevrenderpams}. In the case of what we
call ''spherical sampling'' our approach provides a formula for the recovery
of a function $f$ from a ''band-limited class $PV_{0}$'' by taking its values
over the spheres centered at the origin of $\mathbb{R}^{n}$ and having radii
$e^{j}$ for all $j\in\mathbb{Z}.$ In the second case, of what we call
''parallel sampling'' we recover a function $f$ from a ''band-limited class
$PV_{0}$'' by taking its values on all hyperplanes in $\mathbb{R}^{n}$ defined
by $\left\{  x_{1}=j\right\}  $ for $j\in\mathbb{Z}.$

An interesting feature of our results is that they use essentially all
advances in the one-dimensional sampling theorems for Riesz basis as developed
by Gilbert Walter in his book \cite{walter}, see also the more recent paper
\cite{walterpaper}.

\section{Shannon-Walter sampling with exponential splines ($L-$splines)}

In the present Section we will construct a formula of Shannon type which is
based on exponential splines (these splines are called sometimes $L-$splines),
cf. \cite{okbook}.

\subsection{Shannon sampling for Riesz basis according to Gilbert Walter}

First, we will provide the construction of Gilbert Walter of sampling and the
corresponding Shannon-Walter functions for semi-orthogonal Wavelet Analysis
generated by a scaling function from a Riesz basis; the proofs follow the
proofs outlined in \cite{walter} for the case of orthogonal scaling functions.

Let us assume that the real-valued function $\phi\left(  t\right)  $ is
continuous and has shifts $\left\{  \phi\left(  t-j\right)  \right\}
_{j\in\mathbb{Z}}$ which represent a Riesz basis for the Hilbert subspace
$V_{0}$ in $L_{2}\left(  \mathbb{R}\right)  ,$ i.e.
\[
V_{0}:=clos_{L_{2}\left(  \mathbb{R}\right)  }\left\{  \phi\left(  t-j\right)
:j\in\mathbb{Z}\right\}  ,
\]
and there exist two constants $A,B>0$ such that for arbitrary constants
$c_{j}$ holds
\[
A\sum_{j=-\infty}^{\infty}\left|  c_{j}\right|  ^{2}\leq\left\|
\sum_{j=-\infty}^{\infty}c_{j}\phi\left(  t-j\right)  \right\|  _{L_{2}\left(
\mathbb{R}\right)  }^{2}\leq B\sum_{j=-\infty}^{\infty}\left|  c_{j}\right|
^{2};
\]
following the tradition we call $\phi$ scaling function of $V_{0}$.

For every element $f\in V_{0}$ we have
\begin{equation}
f\left(  t\right)  =\sum_{j=-\infty}^{\infty}f_{j}\phi\left(  t-j\right)  .
\label{ft=}%
\end{equation}
Assume that the space $V_{0}$ has the property that for every $f\in V_{0}$
with
\[
f\left(  j\right)  =0\qquad\text{for }j\in\mathbb{Z}%
\]
follows $f\equiv0.$ Then there should exist a basis $\left\{  \psi_{j}\left(
t\right)  \right\}  _{j\in\mathbb{Z}}$ of $V_{0}$ such that a formula of
Shannon type holds, i.e. for every $f\in V_{0}$
\[
f\left(  t\right)  =\sum_{j=-\infty}^{\infty}f\left(  j\right)  \psi
_{j}\left(  t\right)  \qquad\text{for }t\in\mathbb{R}.
\]
The problem of finding the proper conditions on $\phi$ and the basis $\left\{
\psi_{j}\left(  t\right)  \right\}  _{j\in\mathbb{Z}}$ has been resolved by
Gilbert Walter, \cite{walter}. Below we provide in detail his construction.

It is well-known that to the Riesz basis $\left\{  \phi\left(  t-j\right)
\right\}  _{j\in\mathbb{Z}}$ there exists unique dual Riesz basis $\left\{
\widetilde{\phi}\left(  t-j\right)  \right\}  _{j\in\mathbb{Z}}$ where duality
means
\begin{equation}
\left\langle \phi\left(  x-j\right)  ,\widetilde{\phi}\left(  x-\ell\right)
\right\rangle =\delta_{j\ell}, \label{fitilde}%
\end{equation}
cf. \cite{young}, \cite{walter}. Let us assume that the functions $\phi$ and
$\widetilde{\phi}$ satisfy the following asymptotic conditions:%
\begin{align}
\phi\left(  t\right)   &  =O\left(  \left|  t\right|  ^{-1-\varepsilon
}\right)  ,\quad\widetilde{\phi}\left(  t\right)  =O\left(  \left|  t\right|
^{-1-\varepsilon}\right) \label{asympfi}\\
\qquad\text{as }t  &  \longrightarrow\pm\infty,\quad t\in\mathbb{R},\nonumber
\end{align}
or, which is equivalent, there exist positive constants $C_{1}>0$ and
$C_{2}>0$ such that for all $t\in\mathbb{R}$ holds
\begin{align*}
\left|  \phi\left(  t\right)  \right|   &  \leq\frac{C_{1}}{\left(  1+\left|
t\right|  \right)  ^{1+\varepsilon}}\\
\left|  \widetilde{\phi}\left(  t\right)  \right|   &  \leq\frac{C_{2}%
}{\left(  1+\left|  t\right|  \right)  ^{1+\varepsilon}}.
\end{align*}

Then a direct estimate shows that the function
\begin{equation}
q\left(  x,y\right)  =\sum_{j\in\mathbb{Z}}\widetilde{\phi}\left(  x-j\right)
\phi\left(  y-j\right)  \label{qxy}%
\end{equation}
is uniformly convergent on every compact set in $\mathbb{R}^{2}$.

On the other hand for every fixed $x_{1}\in\mathbb{R}$ (or uniformly for
$x_{1}$ from a compact subset in $\mathbb{R}$ ) we have the estimate
\[
\left|  q\left(  x_{1},y\right)  \right|  ^{2}\leq\sum_{j,k\in\mathbb{Z}%
}\left|  \phi\left(  y-j\right)  \phi\left(  y-k\right)  \right|  \frac
{C_{2}^{2}}{\left(  1+\left|  x_{1}-j\right|  \right)  ^{1+\varepsilon}\left(
1+\left|  x_{1}-k\right|  \right)  ^{1+\varepsilon}}.
\]
We have also
\begin{align*}
\int_{\mathbb{R}}\left|  \phi\left(  y-j\right)  \phi\left(  y-k\right)
\right|  dy  &  \leq\int_{\mathbb{R}}\frac{C_{1}}{\left(  1+\left|
y-j\right|  \right)  ^{1+\varepsilon}}\frac{C_{1}}{\left(  1+\left|
y-k\right|  \right)  ^{1+\varepsilon}}dy\\
&  \leq C_{1}^{2}\int_{\mathbb{R}}\frac{1}{\left(  1+\left|  y-j\right|
\right)  ^{1+\varepsilon}}dy\\
&  =C_{1}^{2}\int_{\mathbb{R}}\frac{1}{\left(  1+\left|  y\right|  \right)
^{1+\varepsilon}}dy<\infty.
\end{align*}
The last implies the convergence of the integral
\[
\int_{\mathbb{R}}\left|  q\left(  x_{1},y\right)  \right|  ^{2}dy.
\]
Thus for every element $f\in L_{2}\left(  \mathbb{R}\right)  $ the integrals%
\[
\int_{\mathbb{R}}q\left(  x_{1},y\right)  f\left(  y\right)  dy
\]
make sense. In a similar way follows that $q\left(  x,y_{1}\right)  \in
L_{2}\left(  \mathbb{R}\right)  $ uniformly for $y_{1}$ from a compact subset
in $\mathbb{R}.$

Thus, from (\ref{ft=}) and the representation for the dual basis we obtain the
following reproduction property for every $f\in V_{0},$
\begin{align}
\int q\left(  x,y\right)  f\left(  y\right)  dy  &  =f\left(  x\right)
,\label{biort1}\\
\int q\left(  x,y\right)  f\left(  x\right)  dx  &  =f\left(  y\right)  .
\label{biort2}%
\end{align}
Let us note that by the definition it is not clear whether the function $q$ is symmetric.

Let us remark that from (\ref{qxy}) it is obvious that
\begin{equation}
q\left(  t-j,0\right)  =q\left(  t,j\right)  . \label{qt-j}%
\end{equation}

We summarize the main results of Gilbert Walter's approach from \cite{walter}
in the following Proposition.

\begin{proposition}
\label{PWalter}Let $\phi\left(  t\right)  $ be a scaling function satisfying
(\ref{asympfi}), i.e. $\phi\left(  t\right)  $ and its dual defined by
(\ref{fitilde}) satisfy $\phi\left(  t\right)  ,\widetilde{\phi}\left(
t\right)  =O\left(  \left|  t\right|  ^{-1-\varepsilon}\right)  $ as
$t\longrightarrow\pm\infty$ for $t\in\mathbb{R}.$ Let us define the function
$\phi^{\ast}$ by
\begin{equation}
\phi^{\ast}\left(  \xi\right)  :=\sum_{j\in\mathbb{Z}}\phi\left(  j\right)
e^{-i\xi j}. \label{ffi*}%
\end{equation}
We assume that the function $\phi^{\ast}\left(  \xi\right)  $ satisfies the
\textbf{non-zero condition}
\begin{equation}
\phi^{\ast}\left(  \xi\right)  \neq0\qquad\text{for }\xi\in\mathbb{R}.
\label{nonzero}%
\end{equation}
Then the following hold:

1. The system of functions $\left\{  q\left(  t,j\right)  \right\}
_{j\in\mathbb{Z}}=\left\{  q\left(  t-j,0\right)  \right\}  _{j\in\mathbb{Z}}$
represents a \textbf{Riesz basis} of $V_{0},$ where $q$ is defined in
(\ref{qxy})

2. The unique \textbf{dual Riesz basis} $\left\{  S_{j}\left(  t\right)
\right\}  _{j\in\mathbb{Z}}$ corresponding to $\left\{  q\left(  t-j,0\right)
\right\}  _{j\in\mathbb{Z}}$ satisfies
\[
S_{j}\left(  t\right)  =S_{0}\left(  t-j\right)  .
\]
Hence $\left\{  S_{0}\left(  t-j\right)  \right\}  _{j\in\mathbb{Z}}$ is a
Riesz basis as well.

3. The following \textbf{Shannon type formula} holds for all $f$ in the space
$V_{0}$:
\begin{equation}
f\left(  t\right)  =\sum_{j\in\mathbb{Z}}f\left(  j\right)  S_{0}\left(
t-j\right)  \qquad\text{for }t\in\mathbb{R}. \label{shannon}%
\end{equation}

4. The function $S_{0}\left(  t\right)  $ has a Fourier transform which
satisfies
\begin{equation}
\widehat{S_{0}}\left(  \xi\right)  =\frac{\widehat{\phi}\left(  \xi\right)
}{\phi^{\ast}\left(  \xi\right)  }\qquad\text{for }\xi\in\mathbb{R}.
\label{shat}%
\end{equation}

\end{proposition}

\begin{remark}
We will call the function $S_{0}$ \textbf{Shannon-Walter} function. Formula
(\ref{shannon}) may be considered as the analog to the
Shannon-Kotelnikov-Whittaker formula whereby the condition for
band-limitedness of $f$ (i.e. the compactness of the support of the Fourier
transform $\widehat{f}\left(  \xi\right)  $ ) is replaced by $f\in V_{0}$. In
the context of Wavelet Analysis the frequency variable $\xi$ of Fourier
Analysis is replaced by the index $j$.
\end{remark}

%

%TCIMACRO{\TeXButton{Proof}{\proof}}%
%BeginExpansion
\proof
%EndExpansion
1. We will prove that the linear operator $T$ which is defined by the map
\[
T\left[  \widetilde{\phi}\left(  \cdot-j\right)  \right]  \left(  x\right)
=q\left(  x,j\right)
\]
has a continuous extension to $V_{0}.$ For this purpose we will consider the
space $\widehat{V_{0}}$ of the Fourier transforms of the functions in $V_{0}.$

Let us consider the operator $T_{1}$ defined on the space $\widehat{V_{0}}$ as
the operator of multiplication by $\overline{\phi^{\ast}\left(  \xi\right)  }%
$, i.e.
\[
T_{1}\left[  f\right]  \left(  \xi\right)  :=f\left(  \xi\right)
\cdot\overline{\phi^{\ast}\left(  \xi\right)  }.
\]
Since by (\ref{asympfi}) the function $\phi$ has a fast decay it follows that
$\phi^{\ast}\left(  \xi\right)  $ is a bounded, continuous and periodic
function and satisfies
\[
\left|  \phi^{\ast}\left(  \xi\right)  \right|  \geq c_{0}>0\qquad\text{for
}\xi\in\mathbb{R}.
\]
Hence, the linear operator $T_{1}$ is bounded and its inverse defined by
multiplication with $\left(  \overline{\phi^{\ast}\left(  \xi\right)
}\right)  ^{-1}$ is also bounded.

Let us see that
\[
T=F^{-1}T_{1}F
\]
where $F$ denotes the Fourier transform and $F^{-1}$ its inverse. Indeed, it
is enough to check this on the elements $\widetilde{\phi}\left(  x-j\right)
.$ We have
\begin{align*}
F\left[  \widetilde{\phi}\left(  \cdot-j\right)  \right]  \left(  \xi\right)
&  =e^{-ij\xi}\widehat{\widetilde{\phi}}\left(  \xi\right) \\
T_{1}F\left[  \widetilde{\phi}\left(  \cdot-j\right)  \right]  \left(
\xi\right)   &  =\overline{\phi^{\ast}\left(  \xi\right)  }e^{-ij\xi}%
\widehat{\widetilde{\phi}}\left(  \xi\right) \\
&  =\left(  \sum_{\ell\in\mathbb{Z}}\phi\left(  \ell\right)  e^{i\xi\ell
}\right)  e^{-ij\xi}\widehat{\widetilde{\phi}}\left(  \xi\right) \\
&  =\widehat{\widetilde{\phi}}\left(  \xi\right)  \cdot\sum_{\ell\in
\mathbb{Z}}\phi\left(  \ell\right)  e^{i\xi\left(  \ell-j\right)  },
\end{align*}
and taking $F^{-1}$ shows finally that $F^{-1}T_{1}F\left[  \widetilde{\phi
}\left(  \cdot-j\right)  \right]  \left(  t\right)  =q\left(  t-j,0\right)  .$

So we have proved that $T$ is a bounded invertible operator, and by a basic
result in \cite{young} (p. $30$ )\ it follows that the system $\left\{
q\left(  t,j\right)  \right\}  _{j\in\mathbb{Z}}$ is a Riesz basis of $V_{0}.$

2. By a theorem in \cite{young} there is a unique Riesz basis $S_{j}\left(
t\right)  $ which is biorthogonal to $\left\{  q\left(  t-j,0\right)
\right\}  _{j\in\mathbb{Z}},$ i.e.
\[
\left\langle q\left(  t-j,0\right)  ,S_{k}\left(  t\right)  \right\rangle
=\delta_{jk}.
\]
On the other hand due to a change of variable $\tau=t-k$ we have
\begin{align*}
\left\langle q\left(  t-j,0\right)  ,S_{0}\left(  t-k\right)  \right\rangle
&  =\int q\left(  t-j,0\right)  \cdot S_{0}\left(  t-k\right)  dt\\
&  =\int q\left(  \tau+k-j,0\right)  \cdot S_{0}\left(  \tau\right)  d\tau\\
&  =\delta_{k-j,0}%
\end{align*}
which shows that the family of functions $\left\{  S_{0}\left(  t-j\right)
\right\}  _{j\in\mathbb{Z}}$ is biorthogonal to $\left\{  q\left(
t-j,0\right)  \right\}  _{j\in\mathbb{Z}}.$ Hence, by the uniqueness it
follows $S_{j}\left(  t\right)  =S_{0}\left(  t-j\right)  .$

3. Equation (\ref{shannon}) is obtained by the expansion of $f\left(
t\right)  $ in the basis $\left\{  S_{0}\left(  t-j\right)  \right\}
_{j\in\mathbb{Z}}.$ Indeed, from
\[
f\left(  t\right)  =\sum_{j}f_{j}S_{0}\left(  t-j\right)
\]
and by the biorthogonality relation (\ref{biort2}) and (\ref{qt-j}) follows
\[
f\left(  j\right)  =\left\langle f,q\left(  t-j,0\right)  \right\rangle
=f_{j}.
\]

4. For proving (\ref{shat}) we need to apply the Shannon expansion
(\ref{shannon}) to the function $\phi$ and to take the Fourier transform.%

%TCIMACRO{\TeXButton{End Proof}{\endproof}}%
%BeginExpansion
\endproof
%EndExpansion

\subsection{The exponential splines ($L-$splines)}

\subsubsection{Preliminaries on exponential splines}

At the beginning of Wavelet Analysis the splines played and important role
(see some history in \cite{daubechies} and \cite{mallat}). In $1991$ Chui and
Wang have constructed the compactly supported spline wavelets, cf.
\cite{chui}. A characteristic feature of this Wavelet Analysis\ is that it is
semi-orthogonal. Further Wavelet Analysis \ using exponential splines has been
initiated in \cite{DeboorDevoreRon} and has been discussed in detail in
\cite{okbook}, \cite{kounchevrenderappr}. We will refer to the monograph
\cite{okbook} for all details of the outline on exponential splines following below.

The exponential splines are defined by means of ordinary differential
operators with constant coefficients given by polynomials
\begin{equation}
L\left(  z\right)  :=\prod_{j=1}^{N}\left(  z-\lambda_{j}\right)  , \label{Lz}%
\end{equation}
where $\lambda_{j}$ are some real constants. For simplicity sake we will
denote by $\Lambda$ the non-ordered vector of the numbers $\lambda_{j}$
\begin{equation}
\Lambda:=\left[  \lambda_{1},\lambda_{2},...,\lambda_{N}\right]  ,
\label{Lambda}%
\end{equation}
where some of the $\lambda_{j}$'s may repeat (the number of repetitions of
$\lambda_{j}$ is the \emph{multiplicity} $\mu_{j}$ of this $\lambda_{j}$ );
this notation is a convenient way to avoid every time explicitly writing the
multiplicities of the $\lambda_{j}$'s. We will have the operator
\begin{equation}
L\left(  \frac{d}{dt}\right)  :=L_{\Lambda}\left(  \frac{d}{dt}\right)
:=\prod_{j=1}^{N}\left(  \frac{d}{dt}-\lambda_{j}\right)  \label{Lddt}%
\end{equation}
and the space $U_{N}$ of analytic solutions defined by
\begin{equation}
U_{N}=U_{N}\left(  \Lambda\right)  :=\left\{  f\left(  t\right)  :L_{\Lambda
}\left(  \frac{d}{dt}\right)  f\left(  t\right)  =0\qquad\text{for }%
t\in\mathbb{R}\right\}  \label{UN}%
\end{equation}
which is well known to be of dimension $N$, cf. \cite{pontryagin}. The space
$U_{N}$ is generated by the exponents
\[
t^{s}e^{\lambda t}\qquad\text{for }s=0,1,...,\mu_{\lambda}-1
\]
where $\ \lambda\in\Lambda$ and $\mu_{\lambda}$ is its multiplicity, where
obviously we have
\[
N=\sum\mu_{\lambda},
\]
and the sum is over all \emph{different} values $\lambda$ in $\Lambda.$

The classical polynomial case is obtained when $\lambda_{j}=0$ for all
$j=1,2,...,N.$

The basic objects which we will consider are the so-called \textbf{cardinal
exponential splines} ($L-$splines ) generated by the operator $L_{\Lambda
}\left(  \frac{d}{dt}\right)  .$ Namely, we define the space $\mathcal{S}%
_{\Lambda}$ by putting
\begin{equation}
\mathcal{S}_{\Lambda}:=\left\{  f:f\left(  t\right)  \in C^{N-2}\left(
\mathbb{R}\right)  \cap C^{\infty}\left(  \mathbb{R}\setminus\mathbb{Z}%
\right)  ,\quad f_{|\left(  j,j+1\right)  }\in U_{N}\right\}  ,
\label{Slambda}%
\end{equation}
i.e. they are piecewise analytic solutions to the equation $L_{\Lambda}\left(
\frac{d}{dt}\right)  f\left(  t\right)  =0.$

We will consider the Wavelet Analysis generated by the $B-$splines arising
from the operator $L_{\Lambda}\left(  \frac{d}{dt}\right)  $ which are called
sometimes $TB-$splines (cf. \cite{schumaker} -- these are cardinal exponential
splines having minimal support). As in the classical polynomial spline theory
there exists up to a factor only one such $TB-$spline, i.e. a cardinal
exponential spline $Q_{N}\left[  \Lambda\right]  \left(  t\right)
=Q_{N}\left(  t\right)  \in\mathcal{S}_{\Lambda}$ supported in the interval
$\left[  0,N\right]  $ and satisfying
\[
Q_{N}\left(  t\right)  >0\qquad\text{for }t\in\left(  0,N\right)  ,
\]
cf. \cite{deboorbook}, \cite{chui}, \cite{schumaker}.

It is convenient to define $Q_{N}\left(  t\right)  $ by means of its Fourier
transform (cf. p. $274$ in \cite{okbook}), namely we have
\begin{equation}
\widehat{Q_{N}}\left(  \xi\right)  =\frac{\prod_{j=1}^{N}\left(
e^{-\lambda_{j}}-e^{-i\xi}\right)  }{\prod_{j=1}^{N}\left(  i\xi-\lambda
_{j}\right)  }\qquad\text{for }\xi\in\mathbb{R}. \label{QNhat}%
\end{equation}

Many important quantities which we will need are defined by means of the
$TB-$spline $Q_{N}\left(  t\right)  .$

In the case of the classical polynomial splines Schoenberg introduced the
so--called Euler-Frobenius polynomial. In order to define a generalization of
the \textbf{Euler-Frobenius polynomials} introduced by \textbf{Schoenberg} we
will consider the function (see Corollary $13.24,$ p. $235$ in \cite{okbook},
Micchelli \cite{Micchelli})
\begin{equation}
A_{N-1}\left(  x;\lambda\right)  :=\frac{1}{2\pi i}\int_{\Gamma}\frac
{1}{L_{\Lambda}\left(  z\right)  }\frac{e^{xz}}{e^{z}-\lambda}dz, \label{AN}%
\end{equation}
where the closed contour $\Gamma$ surrounds all points $\lambda_{j}\in\Lambda$
but excludes all zeros of the function $\frac{e^{xz}}{e^{z}-\lambda}$ . If we
put
\begin{equation}
r\left(  \lambda\right)  :=\prod_{j=1}^{N}\left(  e^{\lambda_{j}}%
-\lambda\right)  ,\qquad s\left(  \lambda\right)  :=\prod_{j=1}^{N}\left(
e^{-\lambda_{j}}-\lambda\right)  \label{rs}%
\end{equation}
then we obtain the following representation (see Corollary $13.25$ p. $235$ in
\cite{okbook})
\begin{equation}
\Pi_{N-1}\left(  x;\lambda\right)  =r\left(  \lambda\right)  A_{N-1}\left(
x;\lambda\right)  \label{PiA}%
\end{equation}
where $\Pi_{N-1}\left(  x;\lambda\right)  $ is a polynomial of degree $\leq
N-1$ of $\lambda.$ The so-called \textbf{Euler-Frobenius} polynomial
\begin{equation}
\Pi_{N-1}\left(  \lambda\right)  :=\Pi_{N-1}\left(  0;\lambda\right)
\label{PiN}%
\end{equation}
has degree $\leq N-2,$ cf. Corollary $13.25$ in \cite{okbook}.

We have the following important Proposition (cf. \cite{Micchelli},
\cite{Schoenberg}, or Theorem $13.31$ on p. $237$ and Corollary $13.53$ on p.
$253$ in \cite{okbook}).

\begin{proposition}
\label{PzerosPiN}The polynomial $\Pi_{N-1}\left(  \lambda\right)  $ has
exactly $N-2$ negative zeros.
\end{proposition}

By means of the $TB-$spline we may define the \textbf{Euler-Schoenberg
exponential spline, }(see p. $254$ in \cite{okbook} ): \
\begin{equation}
\Phi_{N-1}\left(  x;\lambda\right)  =\sum_{j=-\infty}^{\infty}\lambda^{j}%
Q_{N}\left(  x-j\right)  . \label{FiN}%
\end{equation}
Note that the word ''exponential'' has nothing to do with the exponential
meant in the present paper, but it is related to the following easy to check
property
\begin{equation}
\Phi_{N-1}\left(  x+1;\lambda\right)  =\lambda\Phi_{N-1}\left(  x;\lambda
\right)  . \label{Fieq}%
\end{equation}
We have the following important (cf. p. $255,$ Proposition $13.55$ and Theorem
$13.56$ on p. $256,$ Corollary $13.57$ p. $256$ in \cite{okbook}).

\begin{proposition}
\label{PFiPi}1. The Euler-Schoenberg exponential spline in (\ref{FiN})
satisfies the following relation for all $x$ with $0\leq x\leq1,$
\begin{equation}
\Phi_{N-1}\left(  x;\lambda\right)  =\frac{\left(  -1\right)  ^{N-1}}%
{\lambda^{N-1}}e^{-\sum_{j}\lambda_{j}}\Pi_{N-1}\left(  x;\lambda\right)  .
\label{FiPi}%
\end{equation}

2. If the vector $\Lambda$ is symmetric, i.e. satisfies $\Lambda=-\Lambda,$
then
\begin{equation}
\Phi_{N-1}\left(  \frac{N}{2};\frac{1}{z}\right)  =\Phi_{N-1}\left(  \frac
{N}{2};z\right)  \qquad\text{for all }z\in\mathbb{C}, \label{FiNsymmetry}%
\end{equation}

and the polynomial $\Phi_{N-1}\left(  0;z\right)  \neq0$ satisfies
\begin{equation}
\Phi_{N-1}\left(  0;z\right)  \neq0\qquad\text{for all }z\in\mathbb{C}\text{
with }\left|  z\right|  =1. \label{FiNnonzero}%
\end{equation}

\end{proposition}

In particular, we see that from (\ref{Fieq}) and (\ref{FiPi}) the following
representation holds:
\begin{align}
\Phi_{N-1}\left(  \frac{N}{2};\lambda\right)   &  =\lambda^{\frac{N}{2}}%
\Phi_{N-1}\left(  0;\lambda\right) \label{FiPi2}\\
&  =\left(  -1\right)  ^{N-1}\lambda^{-\frac{N}{2}+1}\exp\left(  -\sum
_{j}\lambda_{j}\right)  \Pi_{N-1}\left(  \lambda\right)  .\nonumber
\end{align}

We will use further the following basic

\begin{lemma}
\label{LPiestimate}Let the vector $\Lambda$ be symmetric, i.e. satisfies
$\Lambda=-\Lambda.$ Then the following holds
\[
\left\vert \Pi_{N-1}\left(  -1\right)  \right\vert \leq\left\vert \Pi
_{N-1}\left(  e^{i\xi}\right)  \right\vert \leq\left\vert \Pi_{N-1}\left(
1\right)  \right\vert \qquad\text{for }\xi\in\mathbb{R}.
\]

\end{lemma}

%

%TCIMACRO{\TeXButton{Proof}{\proof}}%
%BeginExpansion
\proof
%EndExpansion
By Proposition \ref{PzerosPiN} the polynomial $\Pi_{N-1}\left(  \lambda
\right)  $ has precisely $N-2$ zeros $v_{j}<0,$ and by Proposition
\ref{PFiPi}, 1) and 2) they satisfy
\[
v_{j}v_{N-1-j}=1\qquad\text{for }j=1,2,...,N-2.
\]
Hence, we have the representation%
\[
\Pi_{N-1}\left(  \lambda\right)  =D\prod_{j=1}^{N-2}\left(  \lambda
-v_{j}\right)  ,
\]
where $D$ is a non-zero coefficient. By the reality of the zeros $v_{j}<0$ it
follows
\begin{align*}
\left|  \Pi_{N-1}\left(  e^{-i\xi}\right)  \right|   &  =\left|  D\right|
\prod_{j=1}^{N-2}\left|  e^{-i\xi}-v_{j}\right|  =\left|  D\right|
\prod_{j=1}^{\frac{N-2}{2}}\left|  e^{-i\xi}-v_{j}\right|  \left|  e^{-i\xi
}-\frac{1}{v_{j}}\right| \\
&  =\left|  D\right|  \prod_{j=1}^{\frac{N-2}{2}}\frac{\left|  e^{-i\xi}%
-v_{j}\right|  ^{2}}{\left|  v_{j}\right|  }=\left|  D\right|  \prod
_{j=1}^{\frac{N-2}{2}}\frac{1-2v_{j}\cos\xi+v_{j}^{2}}{\left|  v_{j}\right|
}\\
&  \leq\left|  D\right|  \prod_{j=1}^{\frac{N-2}{2}}\frac{1-2v_{j}+v_{j}^{2}%
}{\left|  v_{j}\right|  }\\
&  =\left|  \Pi_{N-1}\left(  1\right)  \right|  .
\end{align*}
In a similar way we obtain
\begin{align*}
\left|  \Pi_{N-1}\left(  e^{-i\xi}\right)  \right|   &  \geq\left|  D\right|
\prod_{j=1}^{\frac{N-2}{2}}\frac{1+2v_{j}+v_{j}^{2}}{\left|  v_{j}\right|  }\\
&  =\left|  \Pi_{N-1}\left(  -1\right)  \right|  .
\end{align*}
%

%TCIMACRO{\TeXButton{End Proof}{\endproof}}%
%BeginExpansion
\endproof
%EndExpansion

\subsection{Checking the conditions of Gilbert Walter}

Now we put for the scaling function
\begin{equation}
\phi\left(  t\right)  =Q_{N}\left(  t\right)  . \label{fiQN}%
\end{equation}
A fundamental result in \cite{okbook} (Theorem $14.6$ ) says that for
arbitrary non-ordered vector $\Lambda$ the system $\left\{  \phi\left(
t-j\right)  \right\}  _{j\in\mathbb{Z}}$ represents a Riesz basis for the
space $V_{0}$ defined by means of the cardinal splines $\mathcal{S}_{\Lambda}$
as%
\begin{equation}
V_{0}=\mathcal{S}_{\Lambda}\cap L_{2}\left(  \mathbb{R}\right)  ; \label{V0}%
\end{equation}
the explicit values of the Riesz constants is found in \cite{okbook}.

According to (\ref{QNhat}) the Fourier transform of the basic $TB-$spline
satisfies the asymptotic
\[
\widehat{\phi}\left(  \xi\right)  =\widehat{Q_{N}}\left(  \xi\right)
=\prod_{j=1}^{N}\frac{e^{-\lambda_{j}}-e^{-i\xi}}{i\xi-\lambda_{j}}=O\left(
\frac{1}{\left\vert \xi\right\vert ^{N}}\right)  \qquad\text{for }%
\xi\longrightarrow\infty,
\]
and it is clear that the derivatives of $\widehat{\phi}\left(  \xi\right)  $
satisfy similar asymptotic.

As we mentioned before formula (\ref{fitilde}), the basis $\left\{
\phi\left(  t-j\right)  \right\}  _{j\in\mathbb{Z}}$ has a dual Riesz basis
generated by a unique function $\widetilde{\phi}$, i.e. the set of functions
$\left\{  \widetilde{\phi}\left(  t-j\right)  \right\}  _{j\in\mathbb{Z}}$
generates a Riesz basis of $V_{0}.$ For the dual function $\widetilde{\phi}$
one has the Fourier transform (see p. $356$ in \cite{okbook} )
\begin{align}
\widehat{\widetilde{\phi}}\left(  \xi\right)   &  =C_{S}^{2}\frac
{\widehat{\phi}\left(  \xi\right)  }{\Phi_{2N-1}\left[  \widetilde{\Lambda
}\right]  \left(  N;e^{i\xi}\right)  }\qquad\text{for }\xi\in\mathbb{R}%
\label{fitildehat}\\
\text{with }C_{S}  &  =\exp\left(  \frac{1}{2}\sum_{j=1}^{N}\lambda
_{j}\right)  ,\nonumber
\end{align}
where $\widetilde{\Lambda}=\left[  \Lambda,-\Lambda\right]  $ , i.e. it is the
\textbf{symmetrized vector} of $\Lambda.$

\begin{proposition}
For arbitrary integer $N\geq1$ and arbitrary non-ordered vector $\Lambda$ of
order $N,$ the functions $\phi$ and $\widetilde{\phi}$ defined by (\ref{fiQN})
and (\ref{fitildehat}) satisfy condition (\ref{asympfi}).
\end{proposition}

%

%TCIMACRO{\TeXButton{Proof}{\proof}}%
%BeginExpansion
\proof
%EndExpansion
The asymptotic of the function $\phi$ is clear since $Q_{N}\left(  t\right)  $
is a compactly supported function, in particular, for every integer $m\geq0$
the following asymptotic holds:
\[
\phi\left(  t\right)  =O\left(  \frac{1}{\left\vert t\right\vert ^{m}}\right)
\qquad\text{for }t\longrightarrow\pm\infty.
\]

Now from Lemma \ref{LPiestimate}, applied to the vector $\widetilde{\Lambda}$
and to the corresponding function $\Phi_{2N-1}\left[  \widetilde{\Lambda
}\right]  \left(  N;e^{i\xi}\right)  $ (whereby we use essentially the
symmetry of the vector $\widetilde{\Lambda}$ ), and from formula
(\ref{FiPi2}), it follows that the function $\Phi_{2N-1}\left[  \widetilde
{\Lambda}\right]  \left(  N;e^{i\xi}\right)  $ is bounded, namely
\[
\left|  \Pi_{2N-1}\left(  -1\right)  \right|  \leq\left|  \Phi_{2N-1}\left[
\widetilde{\Lambda}\right]  \left(  N;e^{i\xi}\right)  \right|  \leq\left|
\Pi_{2N-1}\left(  1\right)  \right|  \qquad\text{for }\xi\in\mathbb{R};
\]
here the Euler polynomial $\Pi_{2N-1}\left(  \lambda\right)  $ corresponds to
the vector $\widetilde{\Lambda}.$

On the other hand, we have the following equality proved by integration by
parts:
\begin{align*}
\widetilde{\phi}\left(  t\right)   &  =\frac{1}{2\pi}\int_{-\infty}^{\infty
}e^{i\xi t}\widehat{\widetilde{\phi}}\left(  \xi\right)  d\xi\\
&  =\frac{1}{2\pi}\int_{-\infty}^{\infty}\frac{1}{1+t^{2}}\left(
1-\frac{d^{2}}{d\xi^{2}}\right)  e^{i\xi t}\cdot\widehat{\widetilde{\phi}%
}\left(  \xi\right)  d\xi\\
&  =\frac{1}{2\pi}\frac{1}{1+t^{2}}\int_{-\infty}^{\infty}e^{i\xi t}\left(
1-\frac{d^{2}}{d\xi^{2}}\right)  \widehat{\widetilde{\phi}}\left(  \xi\right)
d\xi.
\end{align*}

Let us estimate the last integral. Since $\Phi_{2N-1}\left[  \widetilde
{\Lambda}\right]  \left(  N;\lambda\right)  $ is a polynomial in $\lambda,$ it
follows that the derivatives of $\Phi_{2N-1}\left[  \widetilde{\Lambda
}\right]  \left(  N;e^{i\xi}\right)  $ with respect to $\xi$ have nice
asymptotic for $\xi\longrightarrow\pm\infty.$ Now by the above remarks about
the asymptotic of the derivatives of the function $\widehat{\phi}\left(
\xi\right)  $ and the representation formula (\ref{fitildehat}) follows the
inequality
\[
\left\vert \widetilde{\phi}\left(  t\right)  \left(  1+t^{2}\right)
\right\vert \leq C\qquad\text{for }t\in\mathbb{R}.
\]
This justifies the asymptotic condition for $\widetilde{\phi}$ in
(\ref{asympfi}).%

%TCIMACRO{\TeXButton{End Proof}{\endproof}}%
%BeginExpansion
\endproof
%EndExpansion

For proving the asymptotic of $\phi$ and $\widetilde{\phi}$ we needed no
restrictions on the vector $\Lambda.$ However the non-zero condition
(\ref{nonzero}) is more restrictive. In fact, G. Walter has shown that the
classical splines (this corresponds to the case of all $\lambda_{j}=0$ in
$\Lambda$ ) generate a \emph{simple and natural} Shannon type formula only in
the case of odd degree cardinal splines, \cite{walter}. Similar is the
situation with the exponential splines. For that reason we will restrict our
attention to the case
\begin{equation}
N=2p.\label{N=2p}%
\end{equation}

\begin{proposition}
\label{Pnonzerosymmetric}

1. We have the equality
\[
\phi^{\ast}\left(  \xi\right)  :=\sum_{j\in\mathbb{Z}}\phi\left(  j\right)
e^{-i\xi j}=\Phi_{N-1}\left[  \Lambda\right]  \left(  0,e^{i\xi}\right)
\qquad\text{for }\xi\in\mathbb{R};
\]
hence the non-zero condition (\ref{nonzero}) is equivalent to the condition
\begin{equation}
\Phi_{N-1}\left[  \Lambda\right]  \left(  0,\lambda\right)  \neq
0\qquad\text{for }\left|  \lambda\right|  =1. \label{nonzeroFi}%
\end{equation}

2. For arbitrary \textbf{symmetric} non-ordered vector $\Lambda$ of order
$N=2p$ and for $\phi$ given by (\ref{fiQN}) condition (\ref{nonzero}) holds,
i.e.
\[
\phi^{\ast}\left(  \xi\right)  \neq0\qquad\text{for }\xi\in\mathbb{R}.
\]

\end{proposition}

%

%TCIMACRO{\TeXButton{Proof}{\proof}}%
%BeginExpansion
\proof
%EndExpansion
By the definition of $\Phi_{N-1}$ in (\ref{FiN}) and by (\ref{FiPi}) follows
\begin{align*}
\sum_{j\in\mathbb{Z}}\phi\left(  j\right)  \lambda^{-j}  &  =\sum
_{j\in\mathbb{Z}}\phi\left(  -j\right)  \lambda^{j}=\Phi_{N-1}\left[
\Lambda\right]  \left(  0,\lambda\right) \\
&  =\frac{\left(  -1\right)  ^{N-1}}{\lambda^{N-1}}\exp\left(  -\sum_{j=1}%
^{N}\lambda_{j}\right)  \Pi_{N-1}\left(  \lambda\right)  .
\end{align*}
(Please note the difference between $\Phi_{N-1}\left[  \Lambda\right]  $ and
$\Phi_{2N-1}\left[  \widetilde{\Lambda}\right]  $ since they correspond to
different vectors). Now Lemma \ref{LPiestimate} gives us the lower and upper
bounds of $\left|  \Pi_{N-1}\left(  e^{i\xi}\right)  \right|  $ from which we
obtain immediately (\ref{nonzero}), i.e.
\[
\sum_{j\in\mathbb{Z}}\phi\left(  j\right)  e^{-i\xi j}\neq0\qquad\text{for
}\xi\in\mathbb{R}.
\]%
%TCIMACRO{\TeXButton{End Proof}{\endproof}}%
%BeginExpansion
\endproof
%EndExpansion

Note that by means of Lemma \ref{LPiestimate} we have used essentially the
symmetry of the vector $\widetilde{\Lambda}$ and the fact that $N=2p.$ We will
apply the above Proposition \ref{Pnonzerosymmetric} to the case of non-ordered
vectors $\Lambda=\Lambda_{k}=\left[  \lambda_{1},...,\lambda_{N}\right]  $
given by
\begin{equation}
\Lambda_{k}:=\left\{
\begin{array}
[c]{c}%
\lambda_{j}=-k\qquad\text{for }j=1,...,p\qquad\\
\lambda_{j}=k\qquad\ \text{for }j=p+1,...,2p
\end{array}
\right.  , \label{Lambda2}%
\end{equation}
where $k$ is some integer.

On the other hand there are other vectors $\Lambda$ which are important for
our multivariate theory, which are not symmetric but close to being symmetric,
and for which the non-zero condition (\ref{nonzero}) holds; they have been
considered in the paper \cite{kounchevrenderJAT} and called \emph{nearly
symmetric} there. For us the following vectors $\Lambda=\Lambda_{k}=\left[
\lambda_{1},...,\lambda_{N}\right]  $ will be important:\
\begin{equation}
\Lambda_{k}:=\left\{
\begin{array}
[c]{c}%
\lambda_{1+j}=k+2j\qquad\text{for }j=0,...,p-1,\qquad\quad\ \ \\
\lambda_{1+p+j}=-n-k+2+2j\qquad\text{for }j=0,...,p-1;
\end{array}
\right.  \label{Lambda1}%
\end{equation}
here $k$ is some integer.

We have the following analog to Proposition \ref{Pnonzerosymmetric}.

\begin{proposition}
\label{PnonzeroNonsymmetric} Let $\Lambda=\Lambda_{k}$ be defined by
(\ref{Lambda1}). Then $\phi^{\ast}$ satisfies the non-zero condition
(\ref{nonzero}).
\end{proposition}

%

%TCIMACRO{\TeXButton{Proof}{\proof}}%
%BeginExpansion
\proof
%EndExpansion
In the proof of Proposition $13$ in the paper \cite{kounchevrenderJAT}, we
proved that
\[
A_{N-1}\left(  0;-1\right)  \neq0
\]
for every $k\geq0.$ On the other hand by Proposition \ref{PzerosPiN} the
polynomial $\Pi_{N-1}\left(  \lambda\right)  $ can have only real zeros with
$\lambda<0$. Since by (\ref{PiA}) and Proposition \ref{PFiPi} we have
\[
\Phi_{N-1}\left[  \Lambda_{k}\right]  \left(  0;\lambda\right)  =-\lambda\cdot
s\left(  \lambda^{-1}\right)  A_{N-1}\left(  0;\lambda\right)
\]
it follows that the function $\Phi_{N-1}\left[  \Lambda_{k}\right]  \left(
0;\lambda\right)  $ has no zeros for $\left|  \lambda\right|  =1.$%

%TCIMACRO{\TeXButton{End Proof}{\endproof}}%
%BeginExpansion
\endproof
%EndExpansion

Hence, in the above cases we may apply the results of G. Walter as presented
in the previous section. By Proposition \ref{PWalter}, 5) it follows that
there exists a \textbf{Shannon-Walter exponential spline} $S_{0}\in
V_{0}=\mathcal{S}_{\Lambda}\cap L_{2}\left(  \mathbb{R}\right)  $ such that
its Fourier transform satisfies
\begin{equation}
\widehat{S_{0}}\left(  \xi\right)  =\frac{\widehat{\phi}\left(  \xi\right)
}{\phi^{\ast}\left(  \xi\right)  }=\frac{\widehat{\phi}\left(  \xi\right)
}{\Phi_{N-1}\left[  \Lambda\right]  \left(  0,e^{i\xi}\right)  }.
\label{S0hat}%
\end{equation}

From the above Propositions we obtain the following:

\begin{theorem}
\label{TShannonExponentialFormula} Let $N=2p$ and $\Lambda$ be symmetric
non-ordered vector of order $N,$ or $\Lambda$ be a vector given by
(\ref{Lambda1}). Then for every cardinal exponential spline $f\in V_{0}$ there
holds a \textbf{Shannon type formula}:
\begin{equation}
f\left(  t\right)  =\sum_{j=-\infty}^{\infty}S_{0}\left(  t-j\right)  f\left(
j\right)  . \label{shannonexpoFormula}%
\end{equation}

\end{theorem}

\section{The multivariate case}

Our main purpose is to find a multivariate generalization of the Shannon-type
formula (\ref{shannon}). For that purpose we will use a Wavelet Analysis which
is genuinely multivariate. We will work in the framework of the polyharmonic
Wavelet Analysis\ introduced and studied in \cite{okbook}. It uses piecewise
polyharmonic functions (polysplines) which generalize the piecewise polynomial
functions (and splines) in the one-dimensional spline theory.

We assume that the integer $p\geq1$ is fixed.

We will consider the solutions to the polyharmonic equation
\[
\Delta^{p}u\left(  x\right)  =0
\]
in annular domains, i.e. in annuli $A_{a,b}=\left\{  x\in\mathbb{R}%
^{n}:a<\left|  x\right|  <b\right\}  $ given by two numbers $a,b>0.$ Here%
\[
\Delta=\sum_{j=1}^{n}\frac{\partial^{2}}{\partial x_{j}^{2}}%
\]
is the Laplace operator and $\Delta^{p}$ is its $p$-th iterate.

For every integer $k\geq0$ we will assume that we are given an orthonormal
basis of the spherical harmonics of homogeneity degree $k,$ which we will
denote as usually by $\left\{  Y_{k,\ell}\left(  x\right)  \right\}  _{\ell
=1}^{a_{k}},$ cf. \cite{steinweiss}, \cite{okbook}. Thus we have
\[
\int_{\mathbb{S}^{n-1}}Y_{k,\ell_{1}}\left(  \theta\right)  Y_{k,\ell_{2}%
}\left(  \theta\right)  d\theta=\delta_{\ell_{1}-\ell_{2}};
\]
here $\mathbb{S}^{n-1}$ is the unit sphere in $\mathbb{R}^{n}$, $d\theta$ is
the spherical area measure, and $\delta_{j}$ is the Kronecker symbol taking on
$1$ for $j=0$ and $0$ elsewhere. We will work with the spherical variables
$x=r\theta$ where $r=\left|  x\right|  ,$ $\theta\in\mathbb{S}^{n-1}.$ It is
known that for special functions of the form $f\left(  r\right)  Y_{k,\ell
}\left(  \theta\right)  $ one has (see \cite{okbook}, p. $152$ )
\[
\Delta^{p}\left(  f\left(  r\right)  Y_{k,\ell}\left(  \theta\right)  \right)
=L_{\left(  k\right)  }^{p}f\left(  r\right)  \cdot Y_{k,\ell}\left(
\theta\right)
\]
where the operator $L_{\left(  k\right)  }$ is given by
\[
L_{\left(  k\right)  }:=\frac{\partial^{2}}{\partial r^{2}}+\frac{n-1}{r}%
\frac{\partial}{\partial r}-\frac{k\left(  n+k-2\right)  }{r^{2}}.
\]
On the other hand by the change $v=\log r$ we obtain another operator
$M_{k,p}$ satisfying
\[
L_{\left(  k\right)  }^{p}\left(  \frac{d}{dr}\right)  =e^{-2pv}M_{k,p}\left(
\frac{d}{dv}\right)
\]
where
\[
M_{k,p}\left(  z\right)  :=\prod_{j=0}^{p-1}\left(  z-k-2j\right)  \left(
z+n+k-2-2j\right)  ,
\]
cf. \cite{okbook}, Theorem $10.34.$ Here it is important that unlike
$L_{\left(  k\right)  }^{p}\left(  \frac{d}{dr}\right)  $ the operator
$M_{k,p}\left(  \frac{d}{dv}\right)  $ has constant coefficients. Thus we are
within the framework of the exponential polynomials and we may consider
exponential splines for the operator $L=M_{k,p}$ defined by the vector
$\Lambda=\Lambda_{k}=\left[  \lambda_{1},...,\lambda_{N}\right]  $ given by
(\ref{Lambda1}).

The following representation for polyharmonic functions in annulus is a
fundamental result (cf. Theorem $10.39$ \ in \cite{okbook}):

\begin{proposition}
\label{PpolyharmonicExpansion} Every function $u\in C^{\infty}\left(
A_{a,b}\right)  $ satisfying $\Delta^{p}u\left(  x\right)  =0$ in the annulus
$A_{a,b}$ permits the representation
\begin{equation}
u\left(  x\right)  =\sum_{k,\ell}u_{k,\ell}\left(  r\right)  Y_{k,\ell}\left(
\theta\right)  ,\qquad\text{for }x=r\theta\label{polyharmonicExpansion}%
\end{equation}
where $u_{k,\ell}\left(  r\right)  $ is a solution to
\[
L_{\left(  k\right)  }^{p}\left(  \frac{d}{dr}\right)  u_{k,\ell}\left(
r\right)  =0\qquad\text{for }a<r<b.
\]

\end{proposition}

Let us assume that the integer $k\geq0$ is fixed.

Further we consider the space of exponential splines generated by the operator
$L=M_{k,p}\left(  \frac{d}{dv}\right)  .$ We define the space of
\emph{cardinal exponential splines} and the Wavelet Analysis\ corresponding to
the operator $L=M_{k,p}\left(  \frac{d}{dv}\right)  $ as in (\ref{V0}), namely
we put
\begin{equation}
\widetilde{V}_{0}^{k}=\left\{  \widetilde{f}\in L_{2}\left(  \mathbb{R}%
\right)  :M_{k,p}\left(  \frac{d}{dv}\right)  \widetilde{f}\left(  v\right)
_{|\left(  j,j+1\right)  }=0,\quad\text{and }\widetilde{f}\in C^{p-2}\left(
\mathbb{R}\right)  \right\}  . \label{V0ktilde}%
\end{equation}

Let us recall that in the case of odd dimension $n$ (since there are no
repeating frequencies in $\Lambda$ of (\ref{Lambda1})) the analytic solutions
of $M_{k,p}\left(  \frac{d}{dv}\right)  f=0$ are linear combinations of the
exponentials
\[
\left\{  e^{\left(  k+2j\right)  v},e^{\left(  -n-k+2+2j\right)  v}\right\}
_{j=0}^{p-1}.
\]

Respectively, by the change of the variable $v=\log r$ we obtain the space
\begin{equation}
V_{0}^{k}:=\left\{  f:f\left(  r\right)  =\widetilde{f}\left(  \log r\right)
,\quad\text{for }\widetilde{f}\in\widetilde{V}_{0}\right\}  . \label{V0k}%
\end{equation}

Finally, we may apply Theorem \ref{TShannonExponentialFormula}, and we will
obtain the \textbf{Shannon-Walter exponential spline} given by (\ref{S0hat})
which we denote further by $\widetilde{S}_{0}^{\left(  k\right)  }\left(
v\right)  .$ By (\ref{shannonexpoFormula}) $\widetilde{S}_{0}^{\left(
k\right)  }\left(  v\right)  $ satisfies a \textbf{Shannon-type formula}, i.e.
for every element $f\in\widetilde{V}_{0}^{\left(  k\right)  }$ we have the
equality
\begin{equation}
\widetilde{f}\left(  v\right)  =\sum_{j=-\infty}^{\infty}\widetilde{f}\left(
j\right)  \widetilde{S}_{0}^{\left(  k\right)  }\left(  v-j\right)  .
\label{shannonexpoKformula}%
\end{equation}
By means of the change $r=e^{v}$ we have $f\left(  e^{v}\right)
:=\widetilde{f}\left(  v\right)  $ and we put
\[
S_{0}^{\left(  k\right)  }\left(  r\right)  :=\widetilde{S}_{0}^{\left(
k\right)  }\left(  \log r\right)  .
\]
Hence, for every $f\in V_{0}^{k}$ we have the \emph{Shannon-type formula}
\begin{equation}
f\left(  r\right)  =\sum_{j=-\infty}^{\infty}f\left(  e^{j}\right)
S_{0}^{\left(  k\right)  }\left(  re^{-j}\right)  .
\label{Shannonexpo-k-formula}%
\end{equation}

\subsection{Asymptotic of $S_{0}^{\left(  k\right)  }$ for $k\longrightarrow
\infty$}

It is essential to see what is the asymptotic of $\widetilde{S}_{0}^{\left(
k\right)  }\left(  v\right)  $ and of $\widehat{\widetilde{S}_{0}^{\left(
k\right)  }}\left(  \xi\right)  $ as $k\longrightarrow\infty.$ This will be
important for the existence and the regularity of the multivariate Shannon
type function.

We assume as above that
\begin{equation}
\phi\left(  t\right)  =Q_{N}\left[  \Lambda\right]  \left(  t\right)
\label{fiQNLambda}%
\end{equation}
for the nonordered vector $\Lambda=\Lambda_{k}=\left[  \lambda_{1}%
,...,\lambda_{N}\right]  $ given by (\ref{Lambda1}) (or (\ref{Lambda2})
below), and that the non-zero condition (\ref{nonzero}) holds.

By formula (\ref{S0hat}) we have
\begin{equation}
\widehat{\widetilde{S}_{0}^{\left(  k\right)  }}\left(  \xi\right)
=\frac{\widehat{\phi_{N}}\left(  \xi\right)  }{\Phi_{N-1}\left[  \Lambda
_{k}\right]  \left(  0;e^{i\xi}\right)  }. \label{S0khat}%
\end{equation}
On the other hand by formula (\ref{FiPi}) we have the representation
\begin{equation}
\Phi_{N-1}\left[  \Lambda_{k}\right]  \left(  0;\lambda\right)  =\frac{\left(
-1\right)  ^{N-1}}{\lambda^{N-1}}\exp\left(  -\sum_{j}\lambda_{j}\right)
\Pi_{N-1}\left(  \lambda\right)  . \label{FiNk}%
\end{equation}

\begin{lemma}
\label{LS0asymptotics} The Fourier transform of the Shannon-Walter function
$S_{0}^{\left(  k\right)  }$ and its Fourier transform satisfy the following asymptotic:\ 

1. uniformly for $\xi\in\mathbb{R}$
\[
\widehat{\widetilde{S}_{0}^{\left(  k\right)  }}\left(  \xi\right)  =O\left(
\frac{1}{k}\right)  \qquad\text{for }k\longrightarrow\infty;
\]

2. uniformly for all $t\in\mathbb{R}$
\[
\widetilde{S}_{0}^{\left(  k\right)  }\left(  t\right)  =O\left(  1\right)
\qquad\text{for }k\longrightarrow\infty.
\]

\end{lemma}

%

%TCIMACRO{\TeXButton{Proof}{\proof}}%
%BeginExpansion
\proof
%EndExpansion
By (\ref{fiQNLambda}) it follows directly that for appropriate constants
$k_{0}\geq0$ and $C>0$ and for all $\xi\in\mathbb{R}$ and $k\geq k_{0}$ we
have
\[
\left|  \widehat{\phi_{N}}\left(  \xi\right)  \right|  =\left|  \frac
{\prod_{j=1}^{N}\left(  e^{-\lambda_{j}}-e^{-i\xi}\right)  }{\prod_{j=1}%
^{N}\left(  i\xi-\lambda_{j}\right)  }\right|  \leq C\frac{e^{pk}}{\left(
\xi^{2}+k^{2}\right)  ^{p}}.
\]

Now we will take care of the denominator in (\ref{S0khat}). We will apply a
subtle estimate of Theorem $11$ in the paper \cite{kounchevrenderJAT} (please
note the conventions in the paper, where the size of the vector $\Lambda$ is
numbered through $N+1$ and here we use $N$ ). It says that there exist two
constants $C_{1},C_{2}>0$ and an integer $k_{0}\geq1$ such that for $k\geq
k_{0}$ holds
\[
\frac{C_{1}}{k^{N-1}}\leq\left\vert A_{N-1}\left(  0;\lambda\right)
\right\vert \leq\frac{C_{2}}{k^{N-1}}\qquad\text{for }\left\vert
\lambda\right\vert =1.
\]
From relation (\ref{PiA}) we obtain
\[
\Pi_{N-1}\left(  \lambda\right)  =r\left(  \lambda\right)  A_{N-1}\left(
0;\lambda\right)
\]
and from (\ref{FiNk}) we obtain
\[
\Phi_{N-1}\left[  \Lambda_{k}\right]  \left(  0;\lambda\right)  =-\lambda\cdot
s\left(  \lambda^{-1}\right)  A_{N-1}\left(  0;\lambda\right)  .
\]
From (\ref{rs}) it follows that the asymptotic of the polynomial $s\left(
\lambda\right)  $ for $k\longrightarrow\infty$ with $\lambda=e^{i\xi}$ is
\[
s\left(  \lambda\right)  \approx e^{kp}.
\]
Hence, for $k\longrightarrow\infty$ and for appropriate constant $C_{3}>0$ we
have
\[
\widehat{\widetilde{S}_{0}^{\left(  k\right)  }}\left(  \xi\right)  \approx
C_{3}\frac{\frac{e^{pk}}{\left(  \xi^{2}+k^{2}\right)  ^{p}}}{\frac{e^{pk}%
}{k^{2p-1}}}=\frac{1}{k}.
\]

2. By definition we have
\[
\widetilde{S}_{0}^{\left(  k\right)  }\left(  t\right)  =\int_{-\infty
}^{\infty}e^{i\xi t}\widehat{S_{0}^{k}}\left(  \xi\right)  d\xi;
\]
hence by a direct estimate and a change of variable $\xi^{\prime}=\xi/k$ we
obtain
\begin{align*}
\left|  \widetilde{S}_{0}^{\left(  k\right)  }\left(  t\right)  \right|   &
\leq C_{3}k^{2p-1}\int_{-\infty}^{\infty}\frac{1}{\left(  \xi^{2}%
+k^{2}\right)  ^{p}}d\xi\\
&  =C_{3}^{\prime}.
\end{align*}
%

%TCIMACRO{\TeXButton{End Proof}{\endproof}}%
%BeginExpansion
\endproof
%EndExpansion

Obviously the results of Lemma \ref{LS0asymptotics} hold for the Shannon
function $S_{0}^{\left(  k\right)  }\left(  r\right)  :=\widetilde{S}%
_{0}^{\left(  k\right)  }\left(  \log r\right)  $ in the variable $r.$

\subsection{The Shannon polyspline}

We have seen that for the vector $\Lambda_{k}$ defined by (\ref{Lambda1}) the
non-zero condition (\ref{nonzeroFi}) is satisfied and we may define the
Shannon-Walter exponential splines $S_{0}^{\left(  k\right)  }.$

We will call the following function (actually it is understood in a
distributional sense) \textbf{Shannon polyspline kernel}:
\begin{align}
S_{0}\left(  r,\theta,\psi\right)   &  :=\sum_{k=0}^{\infty}\sum_{\ell
=1}^{a_{k}}S_{0}^{\left(  k\right)  }\left(  r\right)  Y_{k,\ell}\left(
\theta\right)  Y_{k,\ell}\left(  \psi\right) \label{ShannonPolyspline}\\
&  =\sum_{k=0}^{\infty}S_{0}^{\left(  k\right)  }\left(  r\right)
Z_{k}\left(  \theta,\psi\right) \nonumber\\
&  =\sum_{k=0}^{\infty}\widetilde{S}_{0}^{\left(  k\right)  }\left(  \log
r\right)  Z_{k}\left(  \theta,\psi\right)  ,\nonumber
\end{align}
where $Z_{k}\left(  \theta,\psi\right)  $ is the zonal harmonic of degree $k$,
cf. \cite{steinweiss}. Let us note that for every fixed $\psi$ the function
$S_{0}\left(  r,\theta,\psi\right)  $ is a polyspline of $r\theta$, and for
every fixed $\theta$ the function $S_{0}\left(  r,\theta,\psi\right)  $ is a
polyspline of $r\psi.$

We will characterize the smoothness-convergence properties of this series in
terms of Sobolev spaces. Let us recall that the function $h$ defined by the
series
\[
h\left(  x\right)  :=\sum_{k,\ell}h_{k,\ell}Y_{k,\ell}\left(  \theta\right)
\]
belongs to the Sobolev space on the sphere $H^{s}\left(  \mathbb{S}%
^{n-1}\right)  $ for a real number $s$ if and only if
\[
\sum_{k,\ell}\left(  1+k^{2}\right)  ^{s}h_{k,\ell}^{2}<\infty;
\]
cf. section $7$ and section $22$ in \cite{shubin}.

\begin{lemma}
\label{LShannonSobolev} As a function of $\theta$ (or $\psi$ ) the Shannon
polyspline $S_{0}$ satisfies
\[
S\left(  r,\cdot,\psi\right)  \in H^{s}\left(  \mathbb{S}^{n-1}\right)
,\qquad S\left(  r,\theta,\cdot\right)  \in H^{s}\left(  \mathbb{S}%
^{n-1}\right)
\]
where
\[
s=-n+\frac{3}{2}-\varepsilon
\]
for every number $\varepsilon>0.$
\end{lemma}

%

%TCIMACRO{\TeXButton{Proof}{\proof}}%
%BeginExpansion
\proof
%EndExpansion
We know that for some constant $C>0$ and for all $k\geq0,$ and $\ell
=1,2,...,a_{k}$ holds
\[
\left|  Y_{k,\ell}\left(  \theta\right)  \right|  \leq Ck^{\frac{n}{2}-1},
\]
see e.g. \cite{seeley}, \cite{kalf}. Also, and one has the estimate
\[
a_{k}\leq C_{1}k^{n-2},
\]
cf. e.g. \cite{steinweiss}.

Since the coefficients of the series (\ref{ShannonPolyspline}) are $h_{k,\ell
}=S_{0}^{\left(  k\right)  }\left(  r\right)  Y_{k,\ell}\left(  \psi\right)  $
we have the following inequalities
\[
\sum_{k,\ell}\left(  1+k^{2}\right)  ^{s}h_{k,\ell}^{2}\leq\sum_{k,\ell
}\left(  1+k^{2}\right)  ^{s}C^{2}k^{n-2}\leq C_{2}\sum_{k=0}^{\infty}\left(
1+k^{2}\right)  ^{s}k^{2n-4}%
\]
which proves the theorem.%

%TCIMACRO{\TeXButton{End Proof}{\endproof}}%
%BeginExpansion
\endproof
%EndExpansion

\subsection{The polyharmonic Wavelet Analysis}

In \cite{okbook} we have introduced a \textbf{polyharmonic Wavelet
Analysis}\ by considering the following basic space%
\[
PV_{0}:=clos_{L_{2}\left(  \mathbb{R}^{n}\right)  }\left\{  f\left(  x\right)
:\Delta^{p}f\left(  x\right)  =0\quad\text{on }e^{j}<r<e^{j+1},\quad
x=r\theta,\quad f\in C^{2p-2}\left(  \mathbb{R}^{n}\right)  \right\}  .
\]
Let us note that the radii $e^{j}$ appear in a natural way since by the
representation (\ref{polyharmonicExpansion}) and Proposition
\ref{PpolyharmonicExpansion} every function $u\in PV_{0}$ has coefficients
$u_{k,\ell}\left(  r\right)  $ satisfying
\[
u_{k,\ell}\in V_{0}^{k},
\]
where the spaces $V_{0}^{k}$ are those defined in (\ref{V0k}). Hence, we have
the following decomposition in the sense of $L_{2}$ metric (cf. \cite{okbook}%
)
\[
PV_{0}=\bigoplus_{k=0}^{\infty}\bigoplus_{\ell=1}^{a_{k}}V_{0}^{k}%
\]
which is generated by the expansion (\ref{polyharmonicExpansion}).

As we have seen above, the Shannon polyspline kernel (\ref{ShannonPolyspline})
is a function in some Sobolev space $H^{s}$ and might be a distribution in the
case of a negative exponent. Hence, the formulas which we will write will be
understood in a distributional sense. For two functions (distributions) $f\in
H^{s}\left(  \mathbb{S}^{n-1}\right)  $ and $g\in H^{-s}\left(  \mathbb{S}%
^{n-1}\right)  $ with expansions
\[
f\left(  \theta\right)  =\sum_{k,\ell}f_{k,\ell}Y_{k,\ell}\left(
\theta\right)  ,\qquad g\left(  \theta\right)  =\sum_{k,\ell}g_{k,\ell
}Y_{k,\ell}\left(  \theta\right)  \qquad\text{for }\theta\in\mathbb{S}^{n-1}%
\]
we have the scalar product on $\mathbb{S}^{n-1}$ defined by an integral
understood in a distributional sense
\[
\left\langle f,g\right\rangle :=\int_{\mathbb{S}^{n-1}}f\left(  \theta\right)
g\left(  \theta\right)  d\theta=\sum_{k,\ell}f_{k,\ell}g_{k,\ell}.
\]
Note that the infinite sum which gives sense to the above integral is
absolutely convergent for $f\in H^{s}\left(  \mathbb{S}^{n-1}\right)  $ and
$g\in H^{-s}\left(  \mathbb{S}^{n-1}\right)  .$

In our main Theorem below we will use the spaces $PV_{0}$ to mimic the
one-dimensional cardinal spline space $V_{0},$ and to formulate an analog to
the \textbf{Shannon type formula}.

\begin{theorem}
\label{TShannonpolyspline} For every function $f\in PV_{0}\cap C\left(
\mathbb{R}^{n}\right)  $ the following \textbf{Shannon type formula} holds:
\[
f\left(  r\psi\right)  =\sum_{j=-\infty}^{\infty}\int\limits_{\mathbb{S}%
^{n-1}}S_{0}\left(  \frac{r}{e^{j}},\theta,\psi\right)  f\left(  e^{j}%
\theta\right)  d\theta,
\]
where $S_{0}\left(  r,\theta,\psi\right)  $ is the \textbf{Shannon polyspline
function} defined in (\ref{ShannonPolyspline}), and the integral is understood
in \textbf{distributional sense}.
\end{theorem}

%

%TCIMACRO{\TeXButton{Proof}{\proof}}%
%BeginExpansion
\proof
%EndExpansion
Since by definition $f\in C\left(  \mathbb{R}^{n}\right)  $ it follows that
for all $j\in\mathbb{Z}$ we have $f\left(  e^{j}\theta\right)  \in L_{2},$
hence the following $L_{2}-$convergent expansion for every fixed $r>0$ holds
true
\[
f\left(  x\right)  =\sum_{k=0}^{\infty}\sum_{\ell=1}^{a_{k}}f_{k,\ell}\left(
r\right)  Y_{k,\ell}\left(  \theta\right)  =\sum_{k=0}^{\infty}\sum_{\ell
=1}^{a_{k}}\widetilde{f}_{k,\ell}\left(  \log r\right)  Y_{k,\ell}\left(
\theta\right)  ,
\]
where $x=r\theta.$

We obtain in distributional sense the following equalities, for $v=\log r,$
\begin{align*}
\sum_{j=-\infty}^{\infty}\int\limits_{\mathbb{S}^{n-1}}S_{0}\left(  \frac
{r}{e^{j}},\theta,\psi\right)  f\left(  e^{j}\theta\right)  d\theta &
=\sum_{j=-\infty}^{\infty}\left(  \sum_{k=0}^{\infty}\sum_{\ell=1}^{a_{k}%
}S_{0}^{\left(  k\right)  }\left(  re^{-j}\right)  f_{k\ell}\left(
e^{j}\right)  Y_{k,\ell}\left(  \psi\right)  \right) \\
&  =\sum_{k=0}^{\infty}\sum_{\ell=1}^{a_{k}}\left(  \sum_{j=-\infty}^{\infty
}S_{0}^{\left(  k\right)  }\left(  re^{-j}\right)  f_{k\ell}\left(
e^{j}\right)  Y_{k,\ell}\left(  \psi\right)  \right) \\
&  =\sum_{k=0}^{\infty}\sum_{\ell=1}^{a_{k}}\left(  \sum_{j=-\infty}^{\infty
}\widetilde{S}_{0}^{\left(  k\right)  }\left(  v-j\right)  f_{k\ell}\left(
e^{j}\right)  Y_{k,\ell}\left(  \psi\right)  \right)  .
\end{align*}

On the other hand since for every $k\geq0$ we have $f_{k\ell}\left(
e^{v}\right)  \in V_{0}^{k}$, by the Shannon formula in
(\ref{shannonexpoKformula}) we obtain
\begin{align*}
\sum_{k=0}^{\infty}\sum_{\ell=1}^{a_{k}}\left(  \sum_{j=-\infty}^{\infty
}\widetilde{S}_{0}^{\left(  k\right)  }\left(  v-j\right)  f_{k\ell}\left(
e^{j}\right)  \right)  Y_{k,\ell}\left(  \psi\right)   &  =\sum_{k=0}^{\infty
}\sum_{\ell=1}^{a_{k}}f_{k\ell}\left(  e^{v}\right)  Y_{k,\ell}\left(
\psi\right) \\
&  =f\left(  e^{v}\psi\right) \\
&  =f\left(  r\psi\right)
\end{align*}

This ends the proof.%

%TCIMACRO{\TeXButton{End Proof}{\endproof}}%
%BeginExpansion
\endproof
%EndExpansion

Let us remark that the sampling
\[
\left\{  f\left(  \kappa\right)  :\kappa\in\mathbb{Z}^{n}\right\}
\]
is the one used in the usual multivariate generalizations of the Shannon type
formula. In Theorem \ref{TShannonpolyspline} we use though the sampling
\[
\left\{  f_{k,\ell}\left(  e^{j}\right)  :k\geq0,\ \ell=1,2,...,a_{k}%
,\ j\in\mathbb{Z}\right\}
\]
which defines a new multivariate sampling paradigm.

\section{The polysplines on parallel strips}

We consider very briefly for every integer $p\geq1$ the Wavelet
Analysis\ generated by the polysplines on strips (with periodic conditions),
and the corresponding \emph{Shannon polyspline function} and \emph{Shannon
formula}.

We will introduce some convenient notations:
\begin{align*}
x  &  =\left(  t,y\right)  \in\mathbb{R}^{n}\\
t  &  \in\mathbb{R\quad}\text{and }y\in\mathbb{T}^{n-1},\\
\mathbb{T}^{n-1}  &  :=\left[  0,2\pi\right]  ^{n-1}.
\end{align*}
The \emph{polyspline Wavelet Analysis on strips} is defined by means of the
following space, cf. \cite{okbook},
\[
PV_{0}:=clos_{L_{2}\left(  \mathbb{R}^{n}\right)  }\left\{  f\left(
t,y\right)  :\Delta^{p}f\left(  t,y\right)  =0\quad\text{for }t\notin
\mathbb{Z},\quad f\in C^{2p-2}\left(  \mathbb{R}\times\mathbb{T}^{n-1}\right)
\right\}  ;
\]
these are functions which are piecewise polyharmonic on every strip $\left(
j,j+1\right)  \times\mathbb{T}_{y}^{n-1},$ for all $j\in\mathbb{Z}$ , and
which are $2\pi-$periodic in the variables $y.$ The functions $f\in PV_{0}$
permit the representation
\[
f\left(  t,y\right)  =\sum_{\kappa\in\mathbb{Z}^{n-1}}f_{\kappa}\left(
t\right)  e^{i\left\langle \kappa,y\right\rangle },
\]
where for every $\kappa\in\mathbb{Z}^{n-1}$ the function $f_{\kappa}\left(
t\right)  $ is a solution to the equation
\begin{align*}
\left(  \frac{d^{2}}{dt^{2}}-\left|  \kappa\right|  ^{2}\right)  ^{p}%
f_{\kappa}\left(  t\right)   &  =0\qquad\text{for }t\notin\mathbb{Z},\\
\left|  \kappa\right|  ^{2}  &  =\sum_{j}\kappa_{j}^{2}.
\end{align*}

In the present case we have the vector $\Lambda=\Lambda_{k}=\left[
\lambda_{1},...,\lambda_{N}\right]  $ given by (\ref{Lambda2}) for $k=\left|
\kappa\right|  $. For this $\Lambda=\Lambda_{k}$ and for every integer
$p\geq1$ we obtain the operator
\begin{equation}
L_{\left(  k\right)  }^{p}:=\left(  \frac{d^{2}}{dt^{2}}-k^{2}\right)  ^{p}.
\label{Lkpstrip}%
\end{equation}
Now as in (\ref{V0}) we have similar definition of the cardinal exponential
spline spaces, namely we put
\[
V_{0}^{k}:=\left\{  f\in L_{2}\left(  \mathbb{R}\right)  :L_{\left(  k\right)
}^{p}\left(  \frac{d}{dt}\right)  f\left(  t\right)  _{|\left(  j,j+1\right)
}=0,\quad\text{and }f\in C^{2p-2}\left(  \mathbb{R}\right)  \right\}  .
\]

For the vector $\Lambda_{k}$ defined in (\ref{Lambda2}), by Proposition
\ref{Pnonzerosymmetric} condition (\ref{nonzeroFi}) always holds. Hence, for
every $k\geq0$ there always exists the \textbf{Shannon-Walter exponential
spline} $S_{0}^{\left(  k\right)  }$ defined in (\ref{S0hat}).

\begin{remark}
The case $k=0$ corresponds to the classical cardinal splines of odd degree
considered by G. Walter in \cite{walter}.
\end{remark}

Since all basic lemmata have been proved above also for this type of vectors,
we obtain the following Shannon type formula:

\begin{theorem}
\label{TShannonpolysplineParallel}For every $f\in PV_{0}\cap C\left(
\mathbb{R}^{n}\right)  $ holds the \textbf{Shannon type formula}
\[
f\left(  t,y^{\prime}\right)  =\sum_{j=-\infty}^{\infty}\int_{\mathbb{T}%
^{n-1}}S_{0}\left(  t-j,y,y^{\prime}\right)  f\left(  j,y\right)  dy,
\]
where $S_{0}\left(  t,y,y^{\prime}\right)  $ is the Shannon polyspline
distribution defined by
\[
S_{0}\left(  t,y,y^{\prime}\right)  =\sum_{\kappa\in\mathbb{Z}^{n-1}}%
S_{0}^{\left(  \kappa\right)  }\left(  t\right)  e^{-i\left\langle
\kappa,y\right\rangle }e^{i\left\langle \kappa,y^{\prime}\right\rangle }.
\]

\end{theorem}

Again we may compare the usual multivariate sampling $\left\{  f\left(
\kappa\right)  :\kappa\in\mathbb{Z}^{n}\right\}  $ with our sampling paradigm%
\[
\left\{  f_{\kappa}\left(  j\right)  :\kappa\in\mathbb{Z}^{n-1},j\in
\mathbb{Z}\right\}  .
\]

Finally, an idea for future research might be the error analysis of the
Shannon type formulas of Theorem \ref{TShannonpolyspline} and Theorem
\ref{TShannonpolysplineParallel}.

Corresponding author:

1. Ognyan Kounchev:\ Institute of Mathematics and Informatics, Bulgarian
Academy of Sciences, 8 Acad. G. Bonchev St., 1113 Sofia, Bulgaria

kounchev@math.bas.bg; kounchev@gmx.de

2. Hermann Render: Departamento de Matem\'{a}ticas y Computaci\'{o}n,
Universidad de La Rioja, Edificio Vives, Luis de Ulloa s/n., 26004
Logro\~{n}o, Spain

render@gmx.de

\end{document}